\documentclass{amsart}

\usepackage[colorlinks]{hyperref}

\newtheorem{theorem}{Theorem}[section]

\newtheorem{lemma}{Lemma}[section]

\newtheorem{corollary}{Corollary}[section]

\theoremstyle{remark}
\newtheorem{remark}{Remark}[section]
\theoremstyle{remark}
\newtheorem{remarks}{Remarks}

\theoremstyle{definition}

\newtheorem{definition}{Definition}[section]

\numberwithin{equation}{section}

\renewcommand{\v}[1]{{\mathbf{#1}}}
\newcommand{\inner}[2]{\langle #1,#2\rangle}

\newcommand{\R}{\mathbb{R}}
\newcommand{\mit}{\mathit}
\newcommand{\phhi}{\varphi}
\renewcommand{\phi}{\varphi}

\newcommand{\Kal}[1]{{\mathcal{#1}}}

\def\({\left ( }
\def\){\right )}
\def\<{\left < }
\def\>{\right >}

\begin{document}

\vspace{2cm}

\title[Characterization of the Slant Helix  ...]{Characterization of the Slant Helix as Successor Curve of the General Helix}

\author{Anton Menninger}

\address{Author address: Anton Menninger M.S.\\ Kanbar College of Design, Engineering, and Commerce\\School of Business Administration\newline Philadelphia University\\Philadelphia, PA 19144}
\email{menningera@philau.edu}
\urladdr{\href{http://www.slideshare.net/amenning/presentations}{http://www.slideshare.net/amenning/}}

\subjclass[2010]{53A04}
\keywords{Curve theory; General helix; Slant helix; Curves of constant precession; Frenet equations; Frenet curves; Frenet frame; Bishop frame; Natural equations; Successor curves; Salkowski curves; Closed curves.}

\begin{abstract}
In classical curve theory, the geometry of a curve in three dimensions is essentially characterized by their invariants, curvature and torsion. When they are given, the problem of finding a corresponding curve is known as 'solving natural equations'. Explicit solutions are known only for a handful of curve classes, including notably the plane curves and general helices. 

This paper shows constructively how to solve the natural equations explicitly for an infinite series of curve classes. For every Frenet curve, a family of successor curves can be constructed which have the tangent of the original curve as principal normal. Helices are exactly the successor curves of plane curves and applying the successor transformation to helices leads to slant helices, a class of curves that has received considerable attention in recent years as a natural extension of the concept of general helices. 

The present paper gives for the first time a generic characterization of the slant helix in three-dimensional Euclidian space in terms of its curvature and torsion, and derives an explicit arc-length parametrization of its tangent vector. These results expand on and put into perspective earlier work on Salkowski curves and curves of constant precession, both of which are subclasses of the slant helix. 
\end{abstract}

\maketitle

\section{Introduction}
\let\thefootnote\relax\footnote{This work appeared in \href{http://www.iejgeo.com/matder/dosyalar/makale-170/menninger-2014-7-2-10.pdf}{Int. Electron. J. Geom. 7(2) : 84-91}, 2014.}  
In  classical three-dimensional curve theory, the geometry of a curve is essentially characterized by two scalar functions, {\em curvature} $\kappa$ and {\em torsion} $\tau$, which represent the rate of change of the tangent vector and the osculating plane, respectively. Given two continuous functions of one parameter, there is a space curve (unique up to rigid motion) for which the two functions are its curvature and torsion (parametrized by arc-length). The problem of finding this curve is known as {\em solving natural} or {\em intrinsic equations} and requires solving the Frenet (or an equivalent) system of differential equations (see \cite{scofield}; \cite{Hoppe1862}). One actually solves for the unit tangent and obtains the position vector by integration. Explicit solutions are known only for a handful of curve classes, including notably the plane curves (solved by Euler in 1736) and general helices. 

Recall that {\em general helices} or {\em slope lines} are defined by the property that their tangent makes a constant angle with a fixed direction in every point. Their centrode, the momentary axis of motion, is fixed. Similarly, the principal normal vector of a {\em slant helix} makes a constant angle with a fixed direction and its centrode precesses about a fixed axis. The term was only recently coined (\cite{izumya}) but such curves have appeared in the literature much earlier (e. g. \cite{Hoppe1862}; \cite{bilinski1955}: 180; \cite{Hoschek1969}). In particular, Salkowski studied slant helices with constant curvature (\cite{salkowski1909}, \cite{monterde}) and Scofield derived closed-form arc-length parametrizations for {\em curves of constant precession}, slant helices with the speed of precession (and not just the angle) being constant (\cite{scofield}). 
More recently, there has been considerable interest in the slant helix among curve geometers (e.g. \cite{kula2005}, \cite{kula2010}, \cite{ali}, \cite{Camci2013}).

The close relation between slant helix and general helix leads to the concept of the successor transformation of Frenet curves (\cite{menninger}). Given a Frenet moving frame, we can construct a new Frenet frame in which the original tangent vector plays the role of principal normal. It turns out that general helices are the successor curves of plane curves, and slant helices are the successor curves of general helices. Curves of constant precession are the successor curves of circular helices, which in turn are the successor curves of plane circles. 
This transformation could be repeated indefinitely to give rise to new, yet unexplored classes of curves. 
 
In this paper, the successor transformation of Frenet curves is introduced and its general form given. It is then applied to derive a generic characterization of the slant helix in terms of its curvature and torsion and an explicit arc-length parametrization of its unit tangent vector. 

\section{Basic Notions}

This paper considers space curves in oriented three-dimensional euclidian space $\R^3$. Any regular curve can be represented by an arclength (or {\em unit speed}) parametrization $\v{x}(s)$ with $\v{x}'=\v T$ (e.g. \cite{kuehnel}: 9). Throughout this paper, curves will be assumed to be regular and at least of differentiation order $C^2$, the prime denotes differentiation by the arc length parameter $s$ and $\v T$ denotes the unit tangent vector. 

A {\em moving frame} along a regular curve consists of its unit tangent $\v T$ and two normal unit vector fields $\v N_1, \v N_2$ together forming a {\em positively oriented, orthonormal basis} of $\R^3$ attached to every point of the curve. Two moving frames are called {\em equivalent} if their tangential components are identical. 
Differentiation of a moving frame gives rise to a system of differential equations 
\[ \begin{pmatrix}\v{T}\;\\ \v{N_1}\\ \v{N_2} \end{pmatrix}'
=\begin{pmatrix} \; 0&\;k_1&k_2 \\ -k_1&\;0&k_3 \\ -k_2 &-k_3&\;\,0 \end{pmatrix}\cdot
\begin{pmatrix} \v{T}\; \\ \v{N_1}\\ \v{N_2} \end{pmatrix}\]
with continuous coefficient functions $k_1, k_2, k_3$ arranged in a skew symmetric matrix (due to orthonormality).

The special case $k_2 \equiv 0$ is known as the {\em Frenet moving frame} and $k_3\equiv 0$ as the {\em Bishop moving frame}  (\cite{bishop}). 
In each case, the coefficient functions are sufficient to characterize the associated space curve up to a rigid motion. For the Frenet moving frame, the differential equations take the form of the {\em Frenet equations}
\[\v{T}'=\kappa \v{N}, \qquad \v{N}'=-\kappa \v{T} + \tau \v{B},  \qquad \v{B}'=-\tau \v{N}\]
with unit principal normal $\v{N}$, binormal $\v{B}$, (signed) curvature $\kappa=\inner{\v T'}{ \v N}$ and torsion $\tau = \inner{\v N'}{ \v B}$. 
The {\em Darboux vector} \, $\v D = \tau \v T+\kappa \v B$ represents the Frenet frame's angular momentum. Its direction determines the frame's momentary axis of motion (its {\em centrode}) and its length the angular speed, \,$\omega=\; \sqrt{\kappa^2+\tau^2}$.
The tuple $\Kal F = (\v{T}, \v{N}, \v{B}, \kappa, \tau)$ is called a {\em Frenet apparatus} or  {\em Frenet system} and the pair $(\kappa, \tau)$ a {\em Frenet development} associated with the curve. Further, a Frenet apparatus or moving frame is called {\em generic} if it is stationary along any line segment (i. e. on any interval with vanishing curvature, the torsion also vanishes). 
 
A curve is called a {\em Frenet curve} if it possesses a Frenet frame, and a {\em generic Frenet curve} if it has a generic Frenet frame. This definition differs from convention (e. g. \cite{kuehnel}: 13) in that Frenet curves can have inflection points and even contain line segments and the Frenet frame cannot be assumed unique (cf. \cite{wintner}, \cite{nomizu}, and \cite{wonglai}). However, a generic Frenet apparatus is unique up to the sign of $\v{N},\,\v{B}$, and $\kappa$. Note also that not every $C^2$-curve has a Frenet frame but each has a family of Bishop frames (\cite{bishop}).

\section{The Successor Transformation}

In this section, a transformation of Frenet moving frames is introduced that preserves the Frenet property. 
It is motivated by the observations that (1) any Frenet frame can be transformed into an equivalent Bishop frame, and (2) any Bishop system can be rearranged into a Frenet system. We first provide general formulae for equivalent moving frames, derive the Bishop transformation, and then derive the successor transformation. 

\begin{lemma}[Totality of Moving Frames]
\label{totality}
Let $(\v T, \v N_1, \v N_2)$ be a moving frame of a regular space curve in arclength parametrization. The totality of equivalent moving frames $(\v T,  \v {\overline N_1},  \v {\overline N_2})$ can be expressed as follows, depending on an arbitrary differentiable function $\phhi$:
\[\begin{pmatrix} \v T\; \\ \v {\overline N_1}\\ \v {\overline N_2}\end{pmatrix}=
 \begin{pmatrix} 1&\;0&\;0 \\ 0& \;\cos\phhi& -\sin \phhi  \\ 0&  \;\sin\phhi &  \quad\cos \phhi\end{pmatrix}\cdot
\begin{pmatrix} \v T\; \\ \v {N_1}\\ \v { N_2}\end{pmatrix}\]
The differential equation coefficients get transformed as follows:
\[ \overline k_1=k_1 \cos\phhi - k_2 \sin \phhi, \quad \overline k_2=k_1 \sin\phhi + k_2 \cos \phhi, \quad \overline k_3=k_3-\phhi'.\]

\end{lemma}
\begin{proof}
The rotation matrix, denoted $R_1(\phhi)$, represents the subgroup of the special orthogonal group $SO(3)$ that leaves the first (tangential) component of the moving frame unchanged. It is obvious that all eqivalent tangential moving frames can be expressed in this way. Denoting $F=(\v T, \v N, \v B)^t$ and the coefficient matrix as $K(F)$, we have ${\overline F=R_1(\phhi)F}$, ${F'=K(F)F}$ and \,$\overline F'=R_1'(\phhi)F+R_1(\phhi)F'=
[R_1'(\phhi)+R_1(\phhi)K(F)]R_1(-\phhi)\,\overline F=K(\overline F)\overline F$. Evaluation of $K(\overline F)$ confirms the stated relationships.
\end{proof}

\begin{lemma}[Bishop Transformation]
\label{bishoptransform}
Let $(\v{T}, \v{N}, \v{B}, \kappa, \tau)$ be a Frenet apparatus of a Frenet curve. Then the family of equivalent Bishop apparatuses \,$(\v T, \v N_1, \v N_2, k_1, k_2)$ is given as follows:
\[ \begin{pmatrix} \v {N}_1\\ \v {N}_2\end{pmatrix}=
 \begin{pmatrix} \cos\phhi& -\sin \phhi  \\ \sin\phhi &  \quad\cos \phhi\end{pmatrix}\cdot
\begin{pmatrix} \v {N}\\ \v { B}\end{pmatrix},
\quad 
\begin{pmatrix}  {k}_1\\  { k}_2\end{pmatrix}=\kappa  \begin{pmatrix} \cos\phhi\\ \sin \phhi\end{pmatrix},
\quad \phhi(s)=\phhi_0+\int \tau(s) ds.
\]
depending on a parameter \:$\phhi_0 \in \R$.
This conversion is called {\em Bishop transformation}. 

\end{lemma}
\begin{proof}
The Bishop transformation is an application of lemma \ref{totality} with $k_1=\kappa$, \:$k_2=0$, \:$k_3=\tau$, and $\overline k_3=0$. 
\end{proof}
\newpage
\begin{definition}[Successor Curve]
Given a unit speed $C^2$ curve $\v x: \v x(s)$ with unit tangent $\v T$. A curve $\v x_*:\v x_*(s)$ that has $\v T$ as principal normal is called a {\em successor curve} of $\v x$. 
A Frenet frame $ (\v{T_*}, \v{N_*}, \v{B_*})$ is called  {\em successor frame} of the Frenet frame $(\v{T}, \v{N}, \v{B})$ if $\v N_*\equiv \v T$. 
\end{definition}

\begin{theorem}[Successor Transformation of Frenet Apparatus]
\label{successor}
Every Frenet curve has a family of successor curves. Given a Frenet system $\Kal F=(\v{T}, \v{N}, \v{B}, \kappa, \tau)$, the totality of successor systems  $\Kal F_*=(\v{T}_*, \v{N}_*, \v{B}_*, \kappa_*, \tau_*)$ is as follows:
\[
\begin{pmatrix} \v T_* \\ \v {N_*}\\ \v {B_*}\end{pmatrix}=
 \begin{pmatrix} 0&  -\cos\phhi &  \sin \phhi \\ 1&0&0 \\ 0& \quad\sin\phhi& \cos \phhi  \end{pmatrix}\cdot
\begin{pmatrix} \v T \\ \v {N}\\ \v { B}\end{pmatrix}, \quad
\begin{pmatrix}  \kappa_*\\  \tau_*\end{pmatrix}=\kappa  \begin{pmatrix} \cos\phhi\\ \sin \phhi\end{pmatrix},\quad \phhi(s)\;=\phhi_0+\int \tau(s)ds, 
\]
depending on a parameter \:$\phhi_0 \in \R$. The Darboux vector of $\Kal F_*$ is \,$\v D_*=\kappa \v{B}$.
\end{theorem}

\begin{proof}
Any normal component of a Bishop frame has a tangential derivative and thus is a tangent to a successor curve, and vice versa. Therefore, all successor frames can be obtained by applying the Bishop transformation (lemma \ref{bishoptransform}) to any Frenet frame of the original curve and rearranging the resulting Bishop apparatus \,$(\v T, \v N_1, \v N_2, k_1, k_2)$ into the Frenet apparatus \,$(-\v N_1, \v T, \v N_2, k_1, k_2)$. The family of successor curves is well-defined in that it does not depend on the choice of Frenet system of the original curve. 
\end{proof}

\begin{remark}
The inverse of the successor transformation may be denoted as {\em predecessor transformation.} Bilinski described it for the case $\kappa>0$ (\cite{bilinski1963}) but it is not in general well-defined.
\end{remark}

\begin{corollary}
\label{periodic}
If $\Kal F$ is a periodic Frenet frame of a closed curve, then $\Kal F_*$ is also periodic iff the total torsion of $\Kal F$ is a rational multiple of $\pi$. Note that for a closed curve $c$,  $\int_{c} \tau \mod \pi$ is a curve invariant not depending on the choice of the Frenet frame (\cite{wonglai}: 15).
\end{corollary}

\section{Characterization of the Slant Helix}
In this section, the successor transformation is applied first to a generic plane curve and then to its successor helix to obtain a closed-form expression for the Frenet frame of the slant helix. 

Obviously, a regular $C^2$ space curve lies on a plane if and only if it has a generic Frenet apparatus with constant binormal and vanishing torsion. Given a continuous function $\kappa=\kappa(s)$ defined on an interval $I \subset \R$. Let $\;{\mit\Omega(s)} =\int\kappa(s)ds$ and 
\[\v{T}_P(s)=\begin{pmatrix} \cos{\mit\Omega(s)}\\ \sin{\mit\Omega(s)} \\0 \end{pmatrix}, \quad
\v{N}_P(s)=\begin{pmatrix}   -\sin{\mit\Omega(s)}\\  \quad\cos{\mit\Omega(s)} \\0 \end{pmatrix},\quad
\v{B}_P=\begin{pmatrix} 0\\ 0 \\1 \end{pmatrix}.\]

Then $\Kal F_P= (\v{T}_P, \v{N}_P, \v{B}_P, \kappa, \tau\equiv 0)$ is a generic Frenet apparatus of the plane curve with curvature $\kappa$.

We next turn to the general helix. A proper helix (excluding degenerate cases) makes a constant {\em slope angle}  $\theta\in \; (0, \pi/2)$
with a fixed direction. 
\newpage
\begin{theorem}[Frenet Apparatus of General Helix]\label{helix}
A regular $C^2$ space curve is a general helix if and only if it is a successor curve of a plane curve (and not itself a plane curve). This is equivalent to it having a generic Frenet apparatus satisfying $$\tau_H=\cot\theta \,\kappa_H.$$ 
Given continuous functions ${\kappa_H(s), \tau_H(s)=\cot\theta \, \kappa_H(s)}$. Let ${{\mit\Omega(s)} =\frac{1}{\sin\theta}\int\kappa_H(s)ds}$ and 
\[\v{T}_H(s)=\begin{pmatrix} \quad\sin\theta \sin{\mit\Omega(s)}\\ -\sin\theta\cos{\mit\Omega(s)} \\ \cos\theta \end{pmatrix}, \;
\v{N}_H(s)=\begin{pmatrix}  \cos{\mit\Omega(s)}\\ \sin{\mit\Omega(s)} \\0  \end{pmatrix},\;
\v{B}_H(s)=\begin{pmatrix} -\cos \theta \sin{\mit\Omega(s)}\\ \quad\cos \theta \cos{\mit\Omega(s)}\\ \sin\theta \end{pmatrix}.\]
Then $\Kal F_H= (\v{T}_H, \v{N}_H, \v{B}_H, \kappa_H, \tau_H)$ is a generic Frenet apparatus of a general helix. Its tangent has constant slope $\theta$, its binormal has constant slope $\pi/2-\theta$ and its principal normal image lies on a great circle. Its centrode is fixed and it has angular speed $\kappa$. 
\end{theorem}
\begin{proof}  
Given a helix with unit tangent $\v T$ and $\v D_0$ a unit vector so that $\inner{\v T}{\v D_0}=\cos\theta$. Now $\v T' \perp \v D_0$ and \, $\v N:=\csc\theta\, \v D_0 \times \v T$ is a unit normal vector parallel $\v T'$ and therefore is a principal normal. This construction shows that every $C^2$ helix is a Frenet curve even when inflection points are present. $\int\v N$ is a plane curve with plane normal $\v D_0$ and tangent $\v N$ and by definition, the helix is its successor curve. 

Applying the successor transformation (theorem \ref{successor}) to $\Kal F_P$ with $\phhi(s)=\phhi_0=\pi/2-\theta$ gives $\kappa_H=\sin\theta\,\kappa$ and $\tau_H=\cos\theta\,\kappa$. The Darboux vector is $\kappa \v{B}_P$. $\v T_H$ makes constant angle $\theta$ with $\v{B}_P$ so the successor curve is indeed a helix. Every choice of $\theta\in \; (0, \pi/2)$ results in a distinct proper successor helix. $\Kal F_H$ is generic as $\tau_H$ vanishes whenever $\kappa_H$ vanishes.
\end{proof}

\begin{definition}[Slant Helix]
A Frenet curve is called a {\em slant helix} (\cite{izumya}) if it has a principal normal that has constant slope $\theta\in \; (0, \pi/2)$. 
\end{definition}
Slant helices by definition are exactly the successor curves of proper general helices and a Frenet frame is obtained by applying the successor transformation. A parametrization of the unit tangent is now given. For convenience we'll denote $m:=\cot \theta, \,n:=m/\sqrt{1+m^2}=\cos\theta, \,n/m=\sin\theta$. 

\begin{theorem}[Frenet Apparatus of Slant Helix]\label{slanthelix}
A regular $C^2$ space curve is a slant helix if and only if it has a Frenet development satisfying 
\[\kappa_{SH}=\frac{1}{m} \,\phhi' \cos\phi,\;\quad\tau_{SH}=\frac{1}{m} \,\phhi' \sin\phi\] 
with a differentiable function $\phi(s)$ and $m=\cot\theta\neq 0$.

Given such a Frenet development and let 
$\mit\Omega(s):=\phi(s)/n,\, \lambda_1:=1-n, \lambda_2:= 1+n$ with $n=\cos\theta$. Then the tangent vector of the slant helix thus characterized can be parametrized as follows:
\[\v{T}_{SH}(s)=\frac{1}{2} \begin{pmatrix}
 \lambda_1\cos\lambda_2{\mit\Omega(s)}  +\lambda_2\cos\lambda_1{\mit\Omega(s)}  \\
 \lambda_1\sin\lambda_2{\mit\Omega(s)}  +\lambda_2\sin\lambda_1{\mit\Omega(s)}  \\
    2\frac{n}{m}\sin n{\mit\Omega(s)}   \end{pmatrix}\]

The slant helix is the successor curve of a helix with curvature \,$\kappa_H=\phi'/m$ and torsion \,$\tau_H=\phi'$.
The regular arcs traced by $\v{T}_{SH}$ are spherical helices and the slant helix has a Darboux vector with constant slope $\pi/2-\theta$ and speed of precession $\omega=\kappa_H$. 
Slant helices are generic.
\end{theorem}
\begin{proof}
The successor transformation of $\Kal F_H$ gives \,$\phi=\phi_0+m\int\kappa_H$ and \,$\kappa_{SH}=\kappa_H\cos\phi=m^{-1}\phi'\cos\phi$, $\tau_{SH}=m^{-1}\phi'\sin\phi$. Further, $\mit\Omega=\mit\Omega_0+ m/n\int\kappa_H$. The integration constant $\mit\Omega_0$ (omitted in theorem \ref{helix}) can be set arbitrarily because it only effects a rotation of the frame of the helix. We set $\mit\Omega_0=\phi_0/n$ to get $\phi=n\,\mit\Omega$. 
For the successor tangent, we have $\v{T_{SH}}=-\cos\phhi \v{N}_H + \sin \phhi \v{\v{B}}_H$. Therefore
\[\v{T}_{SH}= \begin{pmatrix} -\cos n\mit\Omega \cos \mit\Omega - n\sin n\mit\Omega\sin\mit\Omega\\
 -\cos n\mit\Omega \sin \mit\Omega + n\sin n \mit\Omega\cos\mit\Omega\\
\frac{n}{m}\sin n \mit\Omega \end{pmatrix}\]
which can be simplified as above (the first two components have been reflected for convenience). Parametrizations of principal normal and binormal can be obtained in the same way but are omitted here. 

The regular segments of the tangent image are slope lines with the principal normal as tangent.
The Darboux vector is $\kappa_H B_H$. The Frenet apparatus thus constructed is generic, since along line segments, $\phhi=const.$ and therefore curvature and torsion both vanish.
\end{proof}

\begin{remarks}
\label{slanthelixrems}
\indent 1. Under the assumption $\kappa \neq 0$, the natural equations of the slant helix are equivalent to the condition
\[ \frac{\kappa^2}{(\kappa^2+\tau^2)^{3/2}} \left( \frac{\tau}{\kappa}\right)' =m.\]
This expression is equivalent to the {\em geodesic curvature} of the spherical image of the principal normal (\cite{izumya}).

2. Solving the differential equation \, $\phhi'\cos\phhi=m\kappa$ for $\phhi$ yields \, $\cos\phhi \,d\phhi=m\,\kappa \,ds\;\Rightarrow \;\sin\phhi(s)=m\int\kappa\, ds$. Similarly, \,$\cos \phi(s)=-m\int\tau ds$. Thus the plane curve with total curvature and total torsion as coordinates is circular. Further, if $\phi \mod 2\pi$ is a periodic function (and consequently the Frenet apparatus is periodic), then total curvature and total torsion over each period vanish. 

3. Restricting \,$\phhi$ to $(-\pi/2, \pi/2)$, the torsion can be expressed as $\tau=\kappa\tan\phhi=\kappa\frac{\sin\phhi}{\sqrt{1-\sin^2\phhi}}$. It follows
\[\tau(s)=\,\kappa(s)\frac{m\int\kappa(s)ds}{\sqrt{1-m^2(\int\kappa(s)ds)^2}}.\]
Given an arbitrary curvature $\kappa$, a torsion function can be constructed to create a slant helix, provided that the domain of $s$ is appropriately restricted. An interesting application is the construction of a {\em Salkowski curve} (\cite{monterde}, \cite{ali}), a slant helix with constant curvature: 
\[\kappa_S(s)\equiv 1, \quad \tau_S(s) = \frac{m s}{\sqrt{1-m^2 s^2}}, \quad m\neq 0, \quad  s\in\, (-1/|m|, +1/|m|).\]
It follows $\,\sin\phi=ms, \,\phi=\arcsin ms, \,\cos\phi=\sqrt{1-m^2}$.
The Salkowski curve is the successor curve of its unit tangent vector, which is a spherical helix and has the same arc length. Its curvatures are
\[\kappa_H=\frac{1}{\sqrt{1-m^2s^2}}, \quad\tau_H=m\,\kappa_H.\]

4. The case \,$\phi'=const.$ characterizes the class of {\em curves of constant precession} (\cite{scofield}) with natural equations
\[\kappa_{CP}=\omega\cos\mu s, \; \tau_{CP} = \omega \sin \mu s, \quad\omega, \mu \neq 0.\]
The curve lies on a hyperboloid of one sheet, and it is closed if and only if $\cos\theta=\mu/\sqrt{\omega^2+\mu^2}$ is rational (figure 2 in \cite{scofield}).

Curves of constant precession are exactly the successor curves of circular helices with curvature $\omega$ and torsion $\mu$, which in turn are successor curves of circles. The Frenet frame of the circular helix has period $L=2\pi\cos\theta/\mu$ and its total torsion per period is $2\pi\cos\theta$. Thus the successor frame is periodic if and only if $\cos\theta$ is rational (corollary \ref{periodic}) and in this case its total curvature and torsion both vanish. 

The successor curve of a closed curve of constant precession again has a periodic Frenet apparatus. It is as yet unknown whether it is closed.
\end{remarks}


\begin{thebibliography}{99}


\bibitem{ali}
Ali, A.~T.
\newblock Position vectors of slant helices in {E}uclidean 3-space.
\newblock {\em Journal of the Egyptian Mathematical Society}, 20:1--6, 2012, DOI: \href{http://www.sciencedirect.com/science/article/pii/S1110256X11000289}{10.1016/j.joems.2011.12.005}.

\bibitem{bilinski1955}
Bilinski, S.
\newblock Eine {V}erallgemeinerung der {F}ormeln von {F}renet und eine {I}somorphie gewisser 
{T}eile der {D}ifferentialgeometrie der {R}aumkurven.
\newblock {\em Glasnik Math.-Fiz. i Astr.}, 10:175--180. 1955

\bibitem{bilinski1963}
Bilinski, S. 
\newblock \"{U}ber eine {E}rweiterungsm\"{o}glichkeit der {K}urventheorie.
\newblock {\em Monatshefte für Mathematik}, 67:289--302, 1963, \href{http://eudml.org/doc/177219}{eudml.org/doc/177219}.

\bibitem{bishop}
Bishop, R.~L.
\newblock There is more than one way to frame a curve.
\newblock {\em Amer. Math. Monthly}, 82:246--251, 1975, \href{http://www.jstor.org/stable/2319846}{www.jstor.org/stable/2319846}.

\bibitem{Camci2013}
Camci, C., Kula, L., and Altinok, M.
\newblock On spherical slant helices in euclidean 3-space.
\href{http://arxiv.org/abs/1308.5532}{arXiv:1308.5532} [math.DG], 2013.

\bibitem{Hoppe1862}
Hoppe, R. 
\newblock \"{U}ber die {D}arstellung der {C}urven durch {K}rümmung und {T}orsion.
\newblock {\em Journal für die reine und angewandte Mathematik}, 60:182--187, 1862, \href{http://eudml.org/doc/147848}{http://eudml.org/doc/147848}.

\bibitem{Hoschek1969}
Hoschek, J. 
\newblock Eine {V}erallgemeinerung der {B}\"oschungsfl\"achen.
\newblock {\em Mathematische Annalen}, 179:275--284, 1969, \href{http://eudml.org/doc/161782}{http://eudml.org/doc/161782}.

\bibitem{izumya}
Izumiya, S. and Takeuchi, N. 
\newblock New {S}pecial {C}urves and {D}evelopable {S}urfaces.
\newblock {\em Turk J Math}, 28:153--163, 2004,
\href {http://journals.tubitak.gov.tr/math/issues/mat-04-28-2/mat-28-2-6-0301-4.pdf}{journals.tubitak.gov.tr/math/issues/mat-04-28-2/mat-28-2-6-0301-4.pdf}.

\bibitem{kuehnel}
K\"uhnel, W. 
\newblock {\em Differential Geometry: Curves - Surfaces - Manifolds, Second Edition}.
\newblock American Mathematical Society, 2006.

\bibitem{kula2010}
{Kula}, L., {Ekmekci}, N., {Yayl{\i}}, Y., and {\.Ilarslan}, K.
\newblock {Characterizations of slant helices in Euclidean 3-space.}
\newblock {\em {Turk. J. Math.}}, 34(2):261--274, 2010,
DOI: \href{http://dx.doi.org/10.3906/mat-0809-17}
{10.3906/mat-0809-17}.

\bibitem{kula2005}
{Kula}, L. and {Yayli}, Y. 
\newblock {On slant helix and its spherical indicatrix.}
\newblock {\em {Appl. Math. Comput.}}, 169(1):600--607, 2010,
DOI: \href{http://dx.doi.org/10.1016/10.1016/j.amc.2004.09.078}
{10.1016/j.amc.2004.09.078}.

\bibitem{menninger}
Menninger, A. 
\newblock {Frenet Curves and Successor Curves.}
\href{http://arxiv.org/abs/1302.3175}{arXiv:1302.3175} [math.DG], 2013.

\bibitem{monterde}
Monterde, J. 
\newblock Salkowski curves revisted: A family of curves with constant curvature
and non-constant torsion.
\newblock {\em Computer Aided Geometric Design}, 26:271–278, 2009, DOI: \href{http://dx.doi.org/10.1016/j.cagd.2008.10.002}
{10.1016/j.cagd.2008.10.002}

\bibitem{nomizu}
Nomizu, K. 
\newblock On {F}renet equations for curves of class ${C}^{\infty}$.
\newblock {\em T\^ohoku Math. J.}, 11:106--112, 1959, \href{http://projecteuclid.org/euclid.tmj/1178244631}{projecteuclid.org/euclid.tmj/1178244631}.

\bibitem{salkowski1909}
Salkowski, E. 
\newblock Zur {T}ransformation von {R}aumkurven. 
\newblock {\em Mathematische Annalen} 
66: 517--557, 1909,
\href{http://eudml.org/doc/158392}{eudml.org/doc/158392}.

\bibitem{scofield}
Scofield, P.~D. 
\newblock \href{http://www.jstor.org/stable/2974768}{Curves of {C}onstant {P}recession.}
\newblock {\em Amer. Math. Monthly}, 102:531--537, 1995.

\bibitem{wintner}
Wintner, A. 
\newblock \href{http://www.jstor.org/stable/2372520}{On {F}renet's {E}quations.}
\newblock {\em Amer. J. Math.}, 78:349--355, 1956.

\bibitem{wonglai}
Wong, Y.-C. and Lai, H.-F. 
\newblock \href{http://projecteuclid.org/euclid.tmj/1178243344}{A {C}ritical {E}xamination} of the {T}heory of {C}urves in {T}hree
  {D}imensional {D}ifferential {G}eometry.
\newblock {\em T\^ohoku Math. J.}, 19:1--31, 1967.


\end{thebibliography}
\end{document}